\documentclass{elsart}
\usepackage{amssymb,amsmath}

\newcommand{\Fin}{\mathrm{Fin}(\beta)}
\newcommand{\db}{\mathrm{d}_{\beta}(1)}
\newcommand{\dbs}{\mathrm{d}_{\beta}^{*}(1)}
\newcommand{\pib}{\pi_{\beta}}
\newcommand{\pia}{\pi_{\alpha}}
\newcommand{\bexp}[1]{\langle #1 \rangle_{\beta}}
\newcommand{\aexp}[1]{\phantom{}_{\alpha}\langle #1 \rangle}
\newcommand{\arep}[1]{\phantom{}_{\alpha}( #1 )}
\newcommand{\lomega}{\phantom{}^{\omega}}
\newcommand{\coloneq}{\mathrel{\mathop:}=}
\newcommand{\eqcolon}{\mathrel={\mathop:}}
\newcommand{\fdot}{\text{\raisebox{0.4ex}[0cm][0cm]{\tiny$\bullet$}}}
\newcommand{\lfdot}{\text{\raisebox{0.2ex}[0cm][0cm]{\tiny$\bullet$}}}

\begin{document}

\begin{frontmatter}
\title{\bf On alpha-adic expansions in Pisot bases\thanksref{words}}

\author[ctu,liafa]{Petr Ambro\v{z}\thanksref{gacr}}
\ead{ampy@linux.fjfi.cvut.cz}
\and
\author[liafa,p8]{Christiane Frougny}
\ead{christiane.frougny@liafa.jussieu.fr}

\address[ctu]{Department of Mathematics, FNSPE, Czech Technical University, Czech Republic}
\address[liafa]{LIAFA, UMR 7089 CNRS \& Universit\'e Paris 7, France}
\address[p8]{Universit\'e Paris 8, France}

\thanks[words]{A preliminary version of the paper was presented at the conference Words 2005~\cite{ambroz-words-05}.}

\thanks[gacr]{Partially supported by Czech Science Foundation GA \v{C}R 201/05/0169 and 
by CTU Internal Grant CTU 0508014}

\begin{abstract}
We study $\alpha$-adic expansions of numbers in an extension field, that is to say, left infinite representations of numbers 
in the positional numeration system with the base $\alpha$, where $\alpha$ is 
an algebraic conjugate of a Pisot number $\beta$.
Based on a result of Bertrand and Schmidt, we prove that a number belongs to $\Qset(\alpha)$
if and only if it has an eventually periodic 
$\alpha$-expansion.
Then we consider $\alpha$-adic expansions of elements of the extension ring $\Zset[\alpha^{-1}]$ 
when $\beta$ satisfies the so-called Finiteness property (F). 
In the particular case that $\beta$ is a quadratic Pisot unit, 
we inspect the unicity and/or multiplicity of $\alpha$-adic expansions of elements of
 $\Zset[\alpha^{-1}]$. We also provide algorithms to generate $\alpha$-adic expansions of rational numbers.
\end{abstract}
\end{frontmatter}
				   
\section{Introduction}
Most usually, real numbers are represented in a 
positional numeration
systems, that is, numbers are considered in the form of finite or
infinite words over a given ordered set --- an alphabet of digits, and their value is taken
following the powers of a real base $\beta>1$. Several different types of
these systems have been studied in the past, e.g.\ usual representations in an
integer base (and its generalizations such as $p$-adic numeration or
systems using signed digits), representations in an irrational base, based
on the so-called $\beta$-expansions (introduced by R\'enyi~\cite{renyi-amash-8}),
or representations with respect to a sequence of integers, like the Fibonacci numeration system.
Another approach is also canonical number systems as studied in~\cite{katai}) for instance.
A survey of most of these concepts was given in~\cite[Chapter 7]{lothaire2}.

In this paper we study another way of the representation of numbers,
strongly connected with the above mentioned representations based on
$\beta$-expansions.
It is called the $\alpha$-adic representation and, roughly speaking,
is a representation of a complex (or real) number in a form of (possibly)
left infinite power series in $\alpha$, where $\alpha$ is a complex (or real) number of modulus less than $1$.

We have two sources of inspiration ---
the $\beta$-numeration systems on one hand and the $p$-adic numbers
(representations of numbers in the form of left infinite power series
in a prime $p$) on the other hand.
However, contrary to the usual $p$-adic numbers the base of the $\alpha$-adic system
is taken to be in modulus smaller than one.
This fact implies an important advantage over the usual $p$-adic
expansions, since we do not have to introduce any special valuation
for the series to converge.

In $\beta$-expansions, numbers are right infinite power series.
The deployment of left infinite power series  
has been used by several authors for different purposes.
Vershik~\cite{vershik-iml-92} (probably the first use of the term
fibadic expansion) and Sidorov and Vershik~\cite{sidorov-vershik-mm-126}
use two-sided expansions to show a connection between symbolic
dynamics of toral automorphisms and arithmetic expansions
associated with their eigenvalues and for study of the
Erd\"os measure (more precisely two-sided generalization of Erd\"os measure).
Two-sided beta-shifts have been studied in full generality by Schmidt in~\cite{schmidt-mm-129}.
Ito and Rao~\cite{ito-rao-pams-04}, and Berth\'e and Siegel~\cite{berthe-siegel-arxiv02}
use representations of two-sided $\beta$-shift in their study of purely periodic expansions
with Pisot unit and non-unit base. The realization by
a finite automaton of the odometer on the two-sided $\beta$-shift has been 
studied by Frougny~\cite{frougny02}.

Left-sided extensions of numeration systems defined by a sequence of
integers, like the Fibonacci numeration system, have been introduced
by Grabner, Liardet and Tichy~\cite{grabner-liardet-tichy-aa-80}, and
studied from the point of view of the odometer function.
The use (at least implicit) of representations
infinite to the left is contained in every study of the Rauzy fractal~\cite{rauzy-bsmf-110},
especially in a study of its border, see e.g.\ Akiyama~\cite{akiyama-dm-2}, Akiyama and
Sadahiro~\cite{akiyama-sadahiro-amiuo-6} or Messaoudi~\cite{messaoudi-aa-95}.

Finally, there is a recent paper by Sadahiro~\cite{sadahiro-pp} on multiply covered points
in the conjugated plane in the case of cubic Pisot units having complex conjugates. Sadahiro's
approach to the left infinite expansions is among all mentioned works the closest one to our
own.

This contribution is organized as follows.
First, we recall known facts about $\beta$-numeration and we define $\alpha$-adic expansions in the case
where $\alpha$ is an algebraic conjugate of a Pisot number $\beta$. 
Recall that, by the results of Bertrand~\cite{bertrand-asp-285}
and Schmidt~\cite{schmidt-blms-12},
a positive real number belongs to the extension field $\Qset(\beta)$ if and only if its
$\beta$-expansion (which is right infinite) is eventually periodic.
Thus it is natural to try to get a similar result for the $\alpha$-adic expansions.
We prove that a number belongs to the field $\Qset(\alpha)$ if and only if
its $\alpha$-adic expansion is eventually periodic to the left with a finite fractional part.
Note that the fields $\Qset(\alpha)$ and $\Qset(\beta)$ are identical, but our result includes also
negative numbers that means one can represent by $\alpha$-adic expansions with positive
digits also negative numbers without utilization of the sign.

Further on, we consider $\alpha$-adic expansions of elements of the ring $\Zset[\alpha^{-1}]$ in the
case when $\beta$ satisfies the Finiteness property (F). We give two algorithms for computing these expansions --- 
one for positive and one for negative numbers. 
Finally, in the case of quadratic Pisot units, we study unicity of the expansions of elements of the ring 
$\Zset[\alpha^{-1}]$. We give 
an algorithm for computing an $\alpha$-adic representation of a rational number and we discuss 
normalization of such representation by means of a finite transducer.

\section{Preliminaries}

\subsection{Words}
An \emph{alphabet} $A$ is a finite ordered set. We denote by the symbol $A^*$ the set of all \emph{finite words} over $A$,
i.e.\ the set of finite concatenation of letters from $A$, empty word (identity of the \emph{free monoid} $A^*$)
is denoted by $\varepsilon$. The set of \emph{infinite words} on $A$ is denoted by $A^{\Nset}$.
A word $u\in A^{\Nset}$ is said to be \emph{eventually periodic} if it is of the form
$u=vz^{\omega}$, where $v,z\in A^*$ are finite words and $z^{\omega} = zzz\cdots$ denotes the infinite concatenation
of $z$ to itself.  We consider also left-infinite words, such a word $u\in \phantom{}^{\Nset}\!A$ is eventually periodic
if $u = \lomega zv$, where $v,z\in A^*$ and $\lomega z  = \cdots zzz$.
A \emph{factor} of a (finite of infinite) word $u$ is a finite word $v$ 
such that $u=v_1vv_2$ for some words $v_1$, $v_2$.

\subsection{Automata and transducers}
An \emph{automaton} over an alphabet $A$, denoted $\mathcal{A} = \langle A,Q,E,I,F \rangle$, is a directed graph 
with labels in $A$. The set $Q$ is set of its vertices, called \emph{states}, $I \subset Q$ is the set 
of \emph{initial states}, $F \subset Q$ is the set of \emph{final states} and 
$E\subset Q\times A\times Q$ is the set of labeled edges, called \emph{transitions}. 
If $(p,a,q)\in E$ one usually writes $p\xrightarrow{a}q$.
The automaton is said to be \emph{finite} if the set of its states is finite.

A \emph{computation} $c$ in $\mathcal{A}$ is a finite sequence of transitions such that
\begin{displaymath}
c = q_0 \xrightarrow{a_1} q_1 \xrightarrow{a_2} q_2\xrightarrow{a_3} \cdots \xrightarrow{a_n} q_n\,.
\end{displaymath}
The \emph{label} of the computation $c$ is a finite word in $A^*$, $a \coloneq a_1a_2\cdots a_n$.
The computation $c$ is \emph{successful} if $q_1 \in I$ and $q_n\in F$. The \emph{behavior} of
$\mathcal{A}$, denoted by $|\mathcal{A}|$, is a subset of $A^*$ of labels of all successful
computations of $\mathcal{A}$. 
An automaton $\mathcal{A}$ is called \emph{deterministic}
if for any pair $(p,a)\in Q\times A$ there exist at most one state $q\in Q$ such that
$p \xrightarrow{a} q$ is a transition of $\mathcal{A}$.

An automaton $\mathcal{T} = \langle A^*\times B^*, Q, E, I, F\rangle$ over a non-free monoid
$A^*\times B^*$ is called a \emph{transducer} from $A^*$ to $B^*$. Its transitions
are labeled by pairs of words $(u,v)\in A^*\times B^*$, the word $u$ is called \emph{input}
and the word $v$ is called \emph{output}. If $(p,(u,v),q)\in E$ one usually writes
$p\xrightarrow{u|v}q$. A computation $c$ in $\mathcal{T}$ is a finite sequence
\begin{displaymath}
c = q_0 \xrightarrow{u_1|v_1} q_1 \xrightarrow{u_2|v_2} q_2\xrightarrow{u_3|v_3} \cdots \xrightarrow{u_n|v_n} q_n\,.
\end{displaymath}
The label of the computation $c$ is $(u,v) \coloneq (u_1u_2\cdots u_n,v_1v_2\cdots v_n)$. 
The behavior of a transducer $\mathcal{T}$ is a relation $R \subset A^*\times B^*$. 
If for any word $u\in A^*$ there exists at most one word $v\in B^*$ such that $(u,v)\in R$ the 
transducer is said to \emph{compute (realize)} a function.
A transducer is called \emph{real-time} if input words of all its transitions are
letters in $A$ (i.e.\ the transitions are labeled in $A\times B^*$).
The \emph{underlying input} (respectively \emph{output}) \emph{automaton}
of a transducer $\mathcal{T}$ is obtained by omitting the output (respectively input) labels
of each transition of $\mathcal{T}$. A transducer is said to be \emph{sequential}
if it is real-time, it has unique initial state and its underlying input automaton
is deterministic. A function is called \emph{sequential} if it can be realized
by a sequential transducer.

\section{Beta expansions}

Let $\beta>1$ be a real number. A \emph{representation in base $\beta$} (or simply a
$\beta$-\emph{representation}) of a real number $x\in\Rset_{+}$ is an infinite sequence
$(x_i)_{i\leq k}$, such that $x_i\in\Zset$ and
\begin{displaymath}
x = x_k\beta^{k} + x_{k-1}\beta^{k-1} + \cdots + x_1\beta + x_0 + x_{-1}\beta^{-1} + x_{-2}\beta^{-2} + \cdots
\end{displaymath}
for a certain $k\in\Zset$.
If a $\beta$-representation of $x$ ends in infinitely many zeros, it is said to be \emph{finite} and the ending zeros
are omitted.

A particular $\beta$-representation --- called \emph{$\beta$-expansion}~\cite{renyi-amash-8} --- is
computed by the so-called greedy algorithm: denote by $\lfloor x \rfloor$, respectively by
$\{x\}$, the integer part, respectively the fractional part, of a number $x$. Find $k\in\Zset$
such that $\beta^{k}\leq x<\beta^{k+1}$. Set $x_k \coloneq \lfloor x/\beta^k \rfloor$
and $r_k\coloneq\{x/\beta^k\}$ and let for $i<k$, $x_i = \lfloor \beta r_{i+1}\rfloor$ and $r_i=\{\beta r_{i+1}\}$.
Then $(x_i)_{i\leq k}$ is the $\beta$-expansion of a number $x$, it is the greatest one
among its $\beta$-representations in the lexicographic order. It is denoted
$\bexp{x}=x_k x_{k-1}\cdots x_0 \lfdot x_{-1} x_{-2} \cdots$,
most significant digit first. When $k$ is negative, we set $x_{-1}=\cdots=x_{-k+1}=0$.
If $\beta$ is not an integer, the digits $x_i$ obtained by the greedy algorithm
are elements of the alphabet $A_{\beta} = \{0, 1, \ldots, \lfloor\beta\rfloor\}$, called the
\emph{canonical alphabet}.

Let $x_k x_{k-1}\cdots x_0 \lfdot x_{-1} x_{-2} \cdots$ be a $\beta$-representation. The \emph{$\beta$-value}
is the function $\pib : A^{\Nset} \mapsto \Rset$ defined by $\pib(x_k x_{k-1}\cdots) \coloneq \sum_{k\geq i} x_i\beta^i$.

Let $C$ be a finite alphabet of integers. The \emph{normalization} on $C$ is the function
$\nu_{C} : C^{\Nset} \rightarrow A_{\beta}^{\Nset}$ that maps a word $w = (w_i)_{i\leq k}$ of $C^{\Nset}$ onto
$\bexp{x}$, where $x=\sum_{i \leq k} w_i\beta^i$, that is, it maps a $\beta$-representation of a number $x$
onto its $\beta$-expansion.

Recall that a \emph{Pisot number} is an algebraic integer $\beta>1$ whose algebraic conjugates are in
modulus less than one.

\begin{thm}[\cite{frougny-mst-25}]
If $\beta$ is a Pisot number then the normalization function $\nu_{C}: C^{\Nset} \rightarrow A_{\beta}^{\Nset}$ 
is computable by a letter-to-letter transducer on any alphabet $C$ of digits.
\end{thm}

A sequence of coefficients which corresponds to some $\beta$-expansion is usually called \emph{admissible}
in the $\beta$-numeration system. For the characterization of admissible sequences we use Parry's
condition~\cite{parry-amash-11}. Let $T_{\beta} : [0,1] \rightarrow [0,1)$ be the 
$\beta$-\emph{transformation} on the unit interval defined by $T_{\beta}(x) \coloneq \{ \beta x \}$. 
The sequence $\db = t_1t_2t_3\cdots$ such that
$t_i = \lfloor \beta T_{\beta}^{i-1}(1)\rfloor$ is called \emph{R\'enyi expansion of 1}. If $\db$ has
infinitely many non-zero digits $t_i$ we set $\dbs = \db$, otherwise if $\ell$ is the greatest
index of non-zero coefficient in $\db$ we set
$\dbs = \left( t_1 t_2 \cdots t_{\ell-1}(t_\ell - 1)\right)^{\omega}$.

\begin{thm}[Parry~\cite{parry-amash-11}]
An infinite sequence $(x_i)_{i\leq k}$ is the $\beta$-expansion of a real number $x\in[0,1)$ if and only if
for all $j\leq k$ the sequence $x_j x_{j-1} x_{j-2} \cdots$ is strictly lexicographically smaller than
the sequence $d_{\beta}^{*}(1)$.
\end{thm}

Properties of $\beta$-expansions are strongly related to symbolic dynamics~\cite{blanchard-tcs-65}.
The closure of the set of admissible $\beta$-expansions is called the $\beta$-shift. It is a symbolic
dynamical system, that is, a closed shift-invariant subset of $\mathcal{A}_{\beta}^{\Nset}$.
A symbolic dynamical system is said to be of \emph{finite type} if the set of its finite factors is
defined by the interdiction of a finite set of words. The $\beta$-shift is of finite type if and only
if $\db$ is finite, see~\cite{blanchard-tcs-65}.

The set of all real numbers $x$ for which the $\beta$-expansion of $|x|$ is finite is denoted
by $\Fin$. 
A number $\beta$ is said to \emph{satisfy Property (F)} if
\begin{displaymath}
\Fin = \Zset[\beta^{-1}]\,.
\end{displaymath}
It has been proved~\cite{frougny-solomyak-etds-12} that Property (F) implies that $\beta$ is a
Pisot number and that $\db$ is  finite.
Conversely, to find a simple algebraic characterization of Pisot numbers satisfying (F) is an
open problem up to now. Let
\begin{equation}\label{eq:min_poly}
\mathrm{M}(\beta) = x^d - a_{d-1} x^{d-1} - \cdots a_1 x - a_0
\end{equation}
be the minimal polynomial of an algebraic integer $\beta$. Several authors have found some 
sufficient conditions on $\mathrm{M}(\beta)$ for $\beta$ to have Property (F).

\begin{thm}[\cite{frougny-solomyak-etds-12}]\label{thm:frosoF}
If the coefficients in~\eqref{eq:min_poly} fulfill
$a_{d-1} \geq a_{d-2} \geq \cdots \geq a_1 \geq a_0 > 0$, then $\beta$ has Property (F).
\end{thm}

\begin{thm}[Hollander~\cite{hollander-thesis}]
If the coefficients in~\eqref{eq:min_poly} fulfill
$a_{d-1} > a_{d-2} + \cdots + a_1 + a_0 > 0$ with $a_i\geq 0$, then $\beta$ has Property (F).
\end{thm}

\begin{thm}[Akiyama~\cite{akiyama-00}]
Let $\beta$ be a cubic Pisot unit. Then $\beta$ has Property (F) if and only if
the coefficients in~\eqref{eq:min_poly} fulfill
$a_0 = 1$, $a_2 \geq 0$ and $-1\leq a_1 \leq a_2+1$.
\end{thm}

Let $\Qset(\beta)$ be the minimal subfield of complex numbers $\Cset$ containing all rationals $\Qset$ as well as
the algebraic number $\beta$. Let $\alpha$ be an algebraic conjugate of $\beta$, then the fields
$\Qset(\beta)$ and $\Qset(\alpha)$ are isomorphic and their isomorphism is induced by the assignment 
$\beta \mapsto \alpha$. Formally, one define isomorphism $': \Qset(\beta) \rightarrow \Qset(\alpha)$
by setting $g(\beta)'=g(\alpha)$, where $g(X)$ is a polynomial in $X$. 

There is a nice characterization of $\beta$-expansions of elements of $\Qset(\beta)$ 
due independently to Bertrand~\cite{bertrand-asp-285} and Schmidt~\cite{schmidt-blms-12}.

\begin{thm}\label{thm:schmidt}
Let $\beta$ be a Pisot number. Then a positive real number $x$ has an eventually periodic $\beta$-expansion if and
only if $x\in\Qset(\beta)$.
\end{thm}

\section{Alpha-adic expansions}

{}From now on let $\beta$ be a Pisot number with finite R\'enyi expansion of 1, say $\db=t_1\cdots t_{\ell}$.
Let $\alpha$ be one of its algebraic conjugates.

\begin{defn}
An \emph{$\alpha$-adic representation} of a number $x\in\Cset$ is a left infinite sequence
$(x_i)_{i\geq -k}$ such that $x_i\in\Zset$ and
\begin{displaymath}
x = \cdots + x_2\alpha^2 + x_{1}\alpha + x_0 + x_{-1}\alpha^{-1} + \cdots + x_{-k}\alpha^{-k}
\end{displaymath}
for a certain $k\in\Zset$. It is denoted $\arep{x} = \cdots x_1 x_0 \lfdot x_{-1} \cdots x_{-k}$.
\end{defn}

\begin{defn}
A (finite, right infinite or left infinite) sequence is said to be \emph{weakly admissible}
if all its finite factors are lexicographically less than or
equal to $\dbs$, which is equivalent to the fact that each factor of length $\ell$ is less than
$t_1\cdots t_{\ell}$ in the lexicographic order. \\
If an $\alpha$-adic representation $(x_i)_{i\geq -k}$ of a number $x$ is
weakly admissible it is said to be an \emph{$\alpha$-adic expansion} of $x$, denoted 
$\aexp{x} = \cdots x_1 x_0 \lfdot x_{-1} \cdots x_{-k}$. 
\end{defn}

\begin{exmp} Let $\beta$ be the golden mean, that is the Pisot number with minimal polynomial
$x^2 - x - 1$. We have $\db = 11$ and $\dbs= (10)^{\omega}$. Hence the sequence 
$(10)^{\omega}$ is a forbidden factor of any $\beta$-expansion. On the other hand,
$\lomega(10)010 \fdot 1$ is an $\alpha$-adic expansion of $-2$.
\end{exmp}

\begin{rem}
Although the $\beta$-expansion of a number is unique, the $\alpha$-adic expansion is not. For instance
in the $\alpha$-adic system associated with the golden mean, the number $-1$ has two $\alpha$-adic expansions
\begin{align*}
\aexp{-1} &= \phantom{1}\lomega(10)\fdot \\
\aexp{-1} &= \lomega(10)0\fdot 1
\end{align*}
\end{rem}

Analogous to the case of $\beta$-representations we define for $\alpha$-adic expansions
the \emph{$\alpha$-value} function $\pia$ and the normalization function $\nu_{C}$.

\section{Eventually periodic $\alpha$-adic expansions}

In order to prove the main theorem about eventually periodic expansions, we need two
technical lemmas.

\begin{lem}\label{lem:pur_per}
Let $y\in(0,1)$ be a real number with the purely periodic $\beta$-expan\-sion
$\bexp{y} = 0\fdot (y_{-1} \cdots y_{-p})^{\omega}$. 
Then $\aexp{-y'} = \lomega(y_{-1} \cdots y_{-p} )\fdot0$.
\end{lem}
\begin{pf}
Suppose
$y= \frac{y_{-1}}{\beta} + \cdots + \frac{y_{-p}}{\beta^p} + \frac{y_{-1}}{\beta^{p+1}} + \cdots$,
which can be also written
$y= \frac{y_{-1}}{\beta} + \cdots + \frac{y_{-p}}{\beta^p} + \frac{y}{\beta^{p}}$. Conjugating the
equation we obtain $y' = \frac{y_{-1}}{\alpha} + \cdots + \frac{y_{-p}}{\alpha^p} + \frac{y'}{\alpha^{p}}$. Hence
$-y' = y_{-1}\alpha^{p-1} + \cdots + y_{-p} - y'\alpha^{p}$ that is $\aexp{-y'}=\lomega(y_{-1}\cdots y_{-p})\fdot\:\!0$,
 which completes the proof.
\qed\end{pf}

\begin{lem}\label{lem:fin}
Let $x\in(0,1)$ be a real number with finite $\beta$-expansion 
$\bexp{x} = 0\fdot x_{-1} \cdots x_{-p}$, then
$\aexp{x'}$ is of the form 
$\lomega(t_1\cdots t_{\ell-1}(t_{\ell}-1)) u_n \cdots u_0 \lfdot u_{-1} \cdots u_{-m}$.
\end{lem}
\begin{pf}
Let $\bexp{x} = 0\fdot x_{-1} \cdots x_{-p}$ with $x_{-p}\neq 0$. By conjugating it and by changing the sign of its
coefficients we obtain an $\alpha$-adic representation of $-x'$, 
$\arep{-x'} = 0\fdot \overline{x_{-1}} \cdots \overline{x_{-p}}$, where $\overline{d}$ denotes the signed
digit $-d$.
If we subtract $-1$ from the last non-zero coefficient $\overline{x_{-p}}$ and replace it by an $\alpha$-adic 
expansion of $-1$, $\aexp{-1} = \lomega(t_1\cdots t_{\ell-1}(t_{\ell}-1))\fdot$, we obtain another
representation, which is eventually periodic with a pre-period of the form of a finite word over the
alphabet $\{-\lfloor\beta\rfloor,\ldots,\lfloor\beta\rfloor\}$. Finally, an $\alpha$-adic expansion
of $-x'$ is simply obtained by the normalization of the pre-period (cf. Algorithm~\ref{algo:z_minus} 
and Example~\ref{ex:z_minus}).
\qed\end{pf}

Lemma~\ref{lem:pur_per} and~\ref{lem:fin} allow us to derive from Theorem~\ref{thm:schmidt} a characterization of
numbers with eventually periodic $\alpha$-adic expansions.
The main difference with Theorem~\ref{thm:schmidt} is that the version for $\alpha$-adic expansions includes also negative 
numbers, that is one can represent by $\alpha$-adic expansions with positive digits also negative numbers without the
necessity of utilization of the sign. 

\begin{thm}
Let $\alpha$ be a conjugate of a Pisot number $\beta$.
A number $x'$ has an eventually periodic $\alpha$-adic expansion
if and only if $x'\in\Qset(\alpha)$.
\end{thm}
\begin{pf}
${\Leftarrow:}$ Let $x'$ have an eventually periodic $\alpha$-adic expansion, say
$\aexp{x'} = \lomega(x_{k+p} \cdots x_{k+1}) x_k \cdots x_0 \lfdot x_{-1} \cdots x_{-j}$.
Let
$u' \coloneq \sum_{i=-j}^k x_i \alpha^i$ and $v' \coloneq \sum_{i={k+1}}^{k+p} x_i \alpha^i$. 
Then $u', v'\in\Zset[\alpha^{-1}]$ and
\begin{displaymath}
x' = u' + \frac {v'}{1 - \alpha^p}\,,
\end{displaymath}
which proves the implication.

${\Rightarrow:}$
Let $x\in \Qset(\beta) \cap [0,1)$.
According to Theorem~\ref{thm:schmidt} the $\beta$-expansion of $x$ is eventually periodic, say
$\bexp{x} = 0\lfdot x_{-1} \cdots x_{-n} ( x_{-n-1} \cdots x_{-n-p} )^{\omega}$. In the case where
the period of $\bexp{x}$ is empty, an eventually periodic $\alpha$-adic expansion of $-x'$ is
obtained by Lemma~\ref{lem:fin}. \\
Let us assume that the period of $\bexp{x}$ is non-empty and let us denote
$y \coloneq \pib(0\fdot ( x_{-(n+1)} \cdots x_{-(n+p)} )^{\omega})$, therefore
$x=\frac{x_{-1}}{\beta} + \cdots \frac{x_{-n}}{\beta^n} + \frac{y}{\beta^n}$.
Conjugating the equation we obtain 
$x'=\frac{x_{-1}}{\alpha} + \cdots + \frac{x_{-n}}{\alpha^n} + \frac{y'}{\alpha^n}$,
hence $-x'= -\frac{y'}{\alpha^n} - \frac{x_{-1}}{\alpha} - \cdots - \frac{x_{-n}}{\alpha^n}$.
According to Lemma~\ref{lem:pur_per} we know how to obtain 
an $\alpha$-adic expansion of $-y'$, hence an $\alpha$-adic
representation of $-x'$ can be obtained by digit wise addition
\begin{displaymath}
\begin{array}{r@{\ }c@{\ }l@{\ }c}
 \lomega(x_{-(n+1)} \cdots x_{-(n+p)} ) x_{-(n+1)} \cdots  x_{-p}\ \lfdot &
x_{-(p+1)} & \cdots & x_{-(n+p)} \\
 \lfdot & (-x_{-1}) & \cdots & (-x_{-n})\\
\hline
 \lomega(x_{-(n+1)} \cdots x_{-(n+p)} ) x_{-(n+1)} \cdots x_{-p}\ \lfdot &
(x_{-(p+1)} -x_{-1} ) & \cdots & ( x_{-(n+p)} - x_{-n})
\end{array}
\end{displaymath}
Therefore we have $\aexp{-x'}$ of the form $\lomega(c_1 \cdots c_p) u$,
where $u$ is a finite word, obtained by the normalization of the pre-period
$ x_{-(n+1)} \cdots x_{-p} \lfdot (x_{-(p+1)} -x_{-1} ) \cdots ( x_{-(n+p)} - x_{-n})$.
Note that this pre-period can be seen as a difference between two finite expansions
and so the normalization will not interfere with the period. \\
Now let $x \ge 1$, $x\in\Qset(\beta)$. Indeed, there exists a positive integer $N$ such
that  $ x < \beta^N$. Hence 
$ t \coloneq 1 - \frac{x}{\beta^N} \in \Qset(\beta) \cap [0,1)$.
As we have proved before the number $ -t = \frac{x'}{\alpha^N} - 1$ has an eventually periodic
$\alpha$-adic expansion. Therefore an eventually periodic $\alpha$-adic expansion $\aexp{x'}$ is simply
obtained by adding $1$ to $\aexp{\frac{x'}{\alpha^N}-1}$, followed by shifting the fractional point
$N$ positions to the left.
\qed \end{pf}

\section{Expansions in bases with Finiteness property (F)}

In the previous section we proved a general theorem characterizing $\alpha$-adic expansions of elements
of the extension field $\Qset(\alpha)$. If we add one additional condition on $\beta$, namely that it
fulfills Property (F), we are able to characterize the expansions of elements of the
ring $\Zset[\alpha^{-1}]$ more precisely.

\begin{prop}\label{prop:z_plus}
Let $\alpha$ be a conjugate of a Pisot number $\beta$ satisfying Property (F).
For any $x\in\Zset[\beta^{-1}]_{+}$ its conjugate $x'$ has at least one $\alpha$-adic expansion.
This expansion is finite and $\aexp{x'}=\bexp{x}$.
\end{prop}
\begin{pf}
Since $\beta$ has Property (F), $\Fin = \Zset[\beta^{-1}]$ and any $x\in\Zset[\beta^{-1}]_{+}$ has a finite 
$\beta$-expansion, say $x = \sum_{i=-j}^{k} x_i\beta^i$. By conjugation we have 
$x' = \sum_{i=-j}^{k} x_i\alpha^i$.
\qed \end{pf}

The proof of Proposition~\ref{prop:z_plus} shows us a way how to compute an $\alpha$-adic
expansion of a number $x'$ which is a conjugate of $x\in\Zset[\beta^{-1}]_{+}$.
The same task is a little bit more complicated in the case where $x'$ is a conjugate of an
$x\in\Zset[\beta^{-1}]_{-}$. An $\alpha$-adic expansion of such a negative number $x'$ is computed by
Algorithm~1 below. 

\begin{alg}\label{algo:z_minus}
Let $x\in\Zset[\beta^{-1}]_{-}$. An $\alpha$-adic expansion of $x'$ is obtained as follows.
\begin{enumerate}
\item Use the greedy algorithm to find the $\beta$-expansion of $-x$, say
$\bexp{-x}  = x_k\cdots x_0\lfdot x_{-1} \cdots x_{-j}$, which is finite
since $\beta$ satisfies Property (F).
\item By changing the signs $x_i \mapsto -x_i$ we obtain an $\alpha$-adic representation
of $x'$ in the form of a finite word over the alphabet $\{0,-1,\ldots,-\lfloor\beta\rfloor\}$.
\item Subtract $-1$ from the rightmost non-zero coefficient $x_{-j}$ and replace it by an
$\alpha$-adic expansion of $-1$, $\aexp{-1} = \lomega(t_1\cdots t_{\ell-1}(t_{\ell}-1))$. 
The representation of $x'$ has now a periodic part $ \lomega(t_1\cdots t_{\ell-1}(t_{\ell}-1))$ 
and a pre-period, which is a finite
word over the alphabet $\{-\lfloor\beta\rfloor,\ldots,\lfloor\beta\rfloor\}$.
\item Finally, the $\alpha$-adic expansion of $x'$ is simply obtained by the normalization of
the pre-period. Note that the pre-period can be seen as a difference between two finite expansions
and so the normalization will not interfere with the period.
\end{enumerate}
\end{alg}

\begin{exmp}\label{ex:z_minus}
Let $\beta$ be the golden mean, $\alpha$ its conjugate. Recall that for example 
$\aexp{-1}=\lomega(10)\fdot$.
We compute an $\alpha$-adic expansion of the number $-4$. The $\beta$-expansion
of $4$ is $101\fdot\, 01$, so $\bar{1}0\bar{1}\fdot\,0\bar{1}$ is an $\alpha$-adic
representation of the number $-4$. Now we subtract $-1$ from the rightmost non-zero coefficient and
replace it by $\aexp{-1}$ as follows
\begin{displaymath}
\begin{array}{r@{}l}
\bar{1}\ 0\ \bar{1} \ \fdot\ & 0\ \bar{1} \\
\ \fdot\ & \phantom{0}\ 1\\
\lomega(1\ 0)\ 1\ 0\ 1\ 0 \ \fdot\ & 1\ 0 \\
\hline
\lomega(1\ 0)\ 1\ \bar{1}\ 1\ \bar{1} \ \fdot\ & 1\ 0
\end{array}
\end{displaymath}
Since the normalization of the pre-period $1 \bar{1} 1 \bar{1}\fdot10 $ gives
$0100\fdot 001$, the expansion is
$\aexp{-4} = \lomega(10)0100\fdot 001$.
\end{exmp}

\begin{prop}\label{prop:z_minus}
Let $\alpha$ be a conjugate of a Pisot number $\beta$ satisfying Property (F).
For any $x\in\Zset[\beta^{-1}]_{-}$, its conjugate $x'$ has at least $\ell$ different $\alpha$-adic expansions, 
which are eventually periodic to the left with the period $\lomega (t_1 \cdots t_{\ell-1}(t_{\ell}-1))$.
\end{prop}
\begin{pf}
First, we show that the number $-1$ has $\ell$ different $\alpha$-adic expansions. 
Recall that $-1 + \pib(\db) = 0$, hence $- \alpha^{\ell} + \alpha^{\ell}\pia(\db) - 1 = - 1$.
Therefore we have the first expansion 
\begin{equation}\label{eq:exp_minus1}
\aexp{-1} = \lomega (t_1 \cdots t_{\ell-1}(t_{\ell}-1))\,.
\end{equation}
Now we successively use the equality $- \alpha^j + \alpha^j\pia(\db) - 1 = -1$ for $j=\ell-1,\ldots, 1$
to obtain  the other $\ell-1$ representations.
For given $j$ this equation is
$-\alpha^j + t_1\alpha^{j-1} + \cdots + t_{j+1}\alpha + (t_j - 1) + t_{j-1}\alpha^{-1} + 
\cdots + t_{\ell}\alpha^{j-\ell} = -1$.
If we replace the coefficient $-1$ at $\alpha^j$ by its expansion~\eqref{eq:exp_minus1} we have 
\begin{equation}\label{eq:exp_minus1_j}
\aexp{-1} = \lomega(t_1 \cdots t_{\ell-1}(t_{\ell}-1))t_1 \cdots t_{j+1} (t_j - 1) \fdot t_{j-1} \cdots t_{\ell}\,.
\end{equation}
Note that periods of expansions obtained in~\eqref{eq:exp_minus1_j} are mutually shifted, they are
situated on all possible $\ell$ positions. That is why all these expansions are essentially distinct.\\
The only difficulty would arise if $t_j = 0$ for some $j$ and hence we would obtain the coefficient $-1$ at $\alpha^0$
by equation~\eqref{eq:exp_minus1_j}. If this is the case we take the pre-period and normalize it
\begin{displaymath}
t_1 \cdots t_{j+1} (t_j - 1) \fdot t_{j-1} \cdots t_{\ell}\ \stackrel{\nu_{C}}{\mapsto}\
u_1 \cdots u_j \fdot u_{j+1} \cdots u_i\,,
\end{displaymath}
where $C = \{-1,0,\ldots,\lfloor\beta\rfloor\}$.

An $\alpha$-adic expansion of $-1$ then will be 
$\aexp{-1} = \lomega(t_1 \cdots t_{\ell-1}(t_{\ell}-1)) u_1 \cdots u_j \lfdot u_{j+1} \cdots u_i$.

Then we consider an $x\in\Zset[\beta^{-1}]_{-}$.
Using the $\ell$ different expansions of $-1$ in Algorithm~\ref{algo:z_minus}
gives us $\ell$ different $\alpha$-adic expansions of the number $x'$. 
\qed\end{pf}
Note that, conversely, if an expansion of a number $z'$ is of the form
$\lomega(t_1 \cdots t_{\ell-1}(t_{\ell}-1)) u \lfdot v$, then $z$ belongs to
$\Zset[\beta^{-1}]_{-}$.

\begin{exmp} Let $\beta$ of minimal polynomial $x^3-x^2-1$;
such a number is Pisot and satisfies the (F) property \cite{akiyama-00}.
We have $\db =101$ and $\dbs= (100)^{\omega}$.
Let $\alpha$ be one of its (complex) conjugates.
The number $-1$ has three different $\alpha$-adic expansion
\begin{align*}
\aexp{-1} &= \phantom{1}\phantom{1}\lomega(100)\fdot \\
\aexp{-1} &= \phantom{1}\lomega(100)0\fdot01 \\
\aexp{-1} &= \lomega(100)01\fdot0000 1
\end{align*}
\end{exmp}

\section{Quadratic Pisot units}

This final section is devoted to quadratic Pisot units, i.e.\ to the algebraic units $\beta$,
with minimal polynomials of the form $x^2-ax-1$, $a\in\Zset_{+}$. Then $\alpha=-\beta^{-1}$.
The R\'enyi expansion of 1 is
$\db=a1$ for such a number $\beta$, which satisfies Property (F), by Theorem~\ref{thm:frosoF}. The canonical
alphabet is $\mathcal{A}=\{0,\ldots, a\}$.

\subsection{Unicity of expansions in $\Zset[\beta]$}
We first establish a technical result.

\begin{prop}\label{lem:fin_num_aexp}
Let $\alpha$ be the conjugate of a quadratic Pisot unit $\beta$.
Let $\phantom{}_{\alpha}\#(x): \Rset \rightarrow \Nset$ be the function counting the number of different
$\alpha$-adic expansions of a number $x$. Then $\phantom{}_{\alpha}\#(x) < +\infty$ for any $x\in\Rset$.
\end{prop}
\begin{pf}
Let $x\in\Rset$ and let $\aexp{x'}=u\lfdot v$ be an $\alpha$-adic expansion of $x'$. Then
\begin{equation}\label{eq:dx}
\begin{split}
& \pib(\fdot v) \in \Zset[\beta] \cap [0,1)\,,\\
& |x' -\pia(\fdot v) |=|\pia(u\fdot )| < \frac{\lfloor\beta\rfloor}{1-|\alpha|}\,.
\end{split}
\end{equation}
Let $D_x \coloneq \{ (\pib(\fdot v),\pia(\fdot v))\ |\ \aexp{x'} = u \lfdot v\}$. Clearly by~\eqref{eq:dx}, 
$D_x$ is a subset of $[0,1) \times \Rset$ with uniformly bounded cardinality,
that is to say there exists a constant $B$ such that $\#D_x \leq B$ for all $x\in\Rset$.

Now suppose that there is a number $y\in\Rset$ such that $y'$ has an infinite number of 
$\alpha$-adic expansions. Indeed, there exists a constant $N$ such that $\alpha^{-N}y'$ 
has $B+1$ different fractional parts. This is in contradiction
with the above proved fact that the number of different fractional parts is uniformly bounded for
$x\in\Rset$.
\qed \end{pf}

Let us note that Proposition~\ref{lem:fin_num_aexp} is conjectured to be valid for all Pisot numbers with
Property (F). In the case that $\beta$ is a cubic Pisot unit with complex conjugates
satisfying Property (F), Sadahiro~\cite{sadahiro-pp} has proved that the above result holds true.

\begin{prop}\label{prop:z_plus2}
Let $\beta$ be a quadratic Pisot unit. Let $x\in\Zset[\beta]_{+}$.
Then $x'$ has a unique $\alpha$-adic expansion, which is finite and such that $\aexp{x'}=\bexp{x}$.
\end{prop}
\begin{pf}
By Proposition~\ref{prop:z_plus} any number $x'\in\Zset[\beta]_{+}$ has an expansion
$\aexp{x'} = x_k \cdots x_0 \lfdot x_{-1} \cdots x_{-j}$.
Let us suppose that $x'$ has another
$\alpha$-adic expansion (necessarily infinite)
$\aexp{x'} = \cdots u_n \cdots u_0 \lfdot u_{-1} \cdots u_{-m}$.
Subtracting these two expansions of $x'$ and normalizing the result we obtain an 
admissible expansion of zero of the form
$\cdots u_{k+3} u_{k+2} v_{k+1}\cdots v_{0} \lfdot v_{-1} \cdots v_{-p}$, with $v_{-p} \neq 0$.
By shifting and relabeling
\begin{equation}\label{eq:exp_zero}
0 = \sum_{i\geq 0} \alpha^i z_i\,,
\end{equation}
where $(z_i)_{i\geq 0}$ is an admissible sequence with $z_0 \neq 0$.
The admissibility condition implies $z_1 \in \{0, \ldots, a-1\}$. Since
$\alpha=-\beta^{-1}$ one can rewrite~\eqref{eq:exp_zero} as
\begin{equation}\label{eq:exp_zero2}
\underbrace{z_0 + \frac{z_2}{\beta^2} + \frac{z_4}{\beta^4} + \cdots}_{\eqcolon LS} = 
\underbrace{\frac{z_1}{\beta} + \frac{z_3}{\beta^3} + \frac{z_5}{\beta^5} + \cdots}_{\eqcolon RS}\,.
\end{equation}
The coefficients $z_i$ for $i\geq 1$ belong to $\{0,\ldots,a\}$, hence
by summing the geometric series on both sides of~\eqref{eq:exp_zero2} we obtain
$LS\in[1,a+\frac{1}{\beta}]$ and $RS\in[0,1-\frac{1}{\beta}]$ which is absurd.
\qed \end{pf}

To prove an analogue of Proposition~\ref{prop:z_plus2} stating the unicity of $\alpha$-adic expansions
for the elements of $\Zset[\beta]_{-}$ we first need the following Lemma.

\begin{lem}
If a number $z$ has an eventually periodic $\alpha$-adic expansion then all its $\alpha$-adic
expansions are eventually periodic.
\end{lem}
\begin{pf}
We have already shown that if a number $x$ has a finite $\alpha$-adic
expansion then this expansion is unique.

Let us consider a number $x'$ with an eventually periodic expansion
\begin{equation}\label{eq:x_per}
\aexp{x'} = \lomega (x_{k+p} \cdots x_{k+1}) x_k \cdots x_0 \lfdot x_{-1} \cdots x_{-j}\,.
\end{equation}
For the sake of contradiction let us assume that $x'$ has another $\alpha$-adic expansion,
which is infinite and non-periodic
\begin{equation}\label{eq:x_non-per}
\aexp{x'} = \cdots u_1 u_0 \lfdot u_{-1} \cdots u_{-m}\,.
\end{equation}
Put $y' \coloneq \alpha^{-(k+1)}x' - \pia(0\lfdot x_k \cdots x_0 x_{-1} \cdots x_{-j})$.
Hence from~\eqref{eq:x_per} we have
\begin{equation}\label{eq:y_per}
\aexp{y'} = \lomega (x_{k+p} \cdots x_{k+1})\fdot 0\,.
\end{equation}
{}From~\eqref{eq:x_non-per}, defining $v_{k+1} \lfdot v_{k} \cdots v_0 v_{-1} \cdots v_{-q}$
as the word
obtained by normalization of the result of digit-wise subtraction
$u_{k+1} \lfdot u_{k} \cdots u_1 u_0 u_{-1} \cdots u_{-m} - 0\lfdot x_k \cdots x_0 x_{-1} \cdots x_{-j}$,
 we have
\begin{equation}\label{eq:y_non-per}
\aexp{y'} = \cdots  u_{k+3} u_{k+2} v_{k+1} \lfdot v_{k} \cdots v_0 v_{-1} \cdots v_{-q}\,,
\end{equation}
which is non-periodic.

Equation~\eqref{eq:y_per} gives us another formula for $y'$,
$y' = \alpha^{-p}y' - \pia(0\lfdot x_{k+p} \cdots x_{k+1})$.
Iterating this formula on the non-periodic expansion~\eqref{eq:y_non-per}
yields infinitely many different $\alpha$-adic expansions of the number $y'$.
This is in the contradiction with the statement of Lemma~\ref{lem:fin_num_aexp}.
\qed \end{pf}

\begin{prop}\label{quad-}
Let $\beta$ be a quadratic Pisot unit.
Let $x\in\Zset[\beta]_{-}$. Then $x'$ has exactly two eventually periodic $\alpha$-adic expansions
with period  $\lomega(a0)$.
\end{prop}
\begin{pf}
At first, we prove that the number $-1$ has no other $\alpha$-adic expansions than those from
Proposition~\ref{prop:z_minus}. Since all $\alpha$-adic expansions of $-1$ have to be eventually
periodic, we will discuss only two cases: when the period is $\lomega(a0)$
and when it is different.
\begin{enumerate}
\item Consider an $\alpha$-adic expansion of $-1$ with the period $\lomega(a0)$
\begin{align*}
\aexp{-1} & = \lomega(a0) d_k \cdots d_0 \lfdot d_{-1} \cdots d_{-j}\,, \\
-1 &= -\alpha^{k+1} + \pia(d_k \cdots d_0 \fdot d_{-1} \cdots d_{-j})\,.
\end{align*} 
The number $-1+\alpha^{k+1}$ is the conjugate of $\beta^{k+1}-1\in\Zset[\beta]_{+}$
and as such has a unique $\alpha$-adic expansion.
Therefore there cannot be two different pre-periods for a given position of the period.
\item Suppose that $-1$ has an $\alpha$-adic expansion with a different period
\begin{displaymath}
\aexp{-1} = \lomega(d_{k+p} \cdots d_{k+1}) d_k \cdots d_0 \lfdot d_{-1} \cdots d_{-j}.
\end{displaymath}
Let $P' \coloneq \pia(d_{k+p} \cdots d_{k+1})$. Then
\begin{displaymath}
-1 = \alpha^{k+1} \frac{P'}{1 - \alpha^p} + \pia(d_k \cdots d_0 \lfdot d_{-1} \cdots d_{-j})\,,
\end{displaymath}
and by taking the conjugate we obtain
\begin{displaymath}
-1 = \beta^{k+1} \frac{P}{1 - \beta^p} + \pib(d_k \cdots d_0 \lfdot d_{-1} \cdots d_{-j})\,.
\end{displaymath}
Therefore
\begin{displaymath}
\underbrace{\pib(d_k \cdots d_0 \lfdot d_{-1} \cdots d_{-j}) + 1}_{\in\Zset[\beta]_{+}}
= \underbrace{\beta^{k+1} \frac{P}{\beta^p - 1}}_{\notin\Zset[\beta]_{+}}\,,
\end{displaymath}
which is a contradiction.
\end{enumerate}

Validity of the statement for numbers $x'\in\Zset[\alpha]_{-}$, $x'\neq -1$, is then a simple consequence
of Algorithm~\ref{algo:z_minus}.
\qed \end{pf}

\subsection{Representations of rational numbers}

In this subsection we inspect $\alpha$-adic expansions of rational numbers. We give below an algorithm
for computing an $\alpha$-adic representation of a rational number $q\in\Qset$, $|q|<1$.
The algorithm for computing $\aexp{q}$
is a sort of a right to left normalization --- it consists of
successive transformations of a representation of $q$, and it gives as a result
a left infinite sequence on the canonical alphabet $\mathcal{A}$. 

Let $x_1$, $x_2$ and $x_3$ be rational digits, and define the following transformation
\begin{equation}\label{eq:alg_trans}
\psi : (x_3)(x_2)(x_1) \mapsto\ (x_3-(\lceil x_1\rceil - x_1))(x_2 +a(\lceil x_1\rceil -x_1))(\lceil x_1\rceil)\,.
\end{equation}
Note that this transformation preserves the $\alpha$-value.

\begin{alg}\label{algo:aexp_q} 

\noindent{\sl Input}:  $q \in \Qset \cap (-1,1)$.\\
{\sl Output}: a sequence  $s=(s_i)_{i \ge 0}$ of $\mathcal{A}^{\Nset}$ such that
$\sum_{i \ge 0} s_i \alpha^{i} =q.$

\noindent{\bf begin}\\
 \hspace*{1cm} $s_0:=q$;\\
 \hspace*{1cm} {\bf for} $i \ge 1$ {\bf do} $s_i:= 0$;\\
 \hspace*{1cm} $i:=0;$\\
  \hspace*{1cm}   {\bf repeat} \\
\hspace*{2cm}  
$s_{i+2} s_{i+1}s_{i}:=\psi(s_{i+2} s_{i+1}s_{i})$; \\
\hspace*{2cm} $i:=i+1;$\\
{\bf end}
\end{alg}

Since 
the starting point of the whole process is a single rational number,
after each step there will be at most two non-integer coefficients --- rational numbers with the
same denominator as $q$.

Denote $s^{(i+1)}$ the resulting sequence after step $i$; thus $s^{(0)}=\lomega0q$ and,
for $i \ge 0$,
$s^{(i+1)}= \cdots s^{(i+1)}_{i+4} s^{(i+1)}_{i+3}s^{(i+1)}_{i+2}s^{(i+1)}_{i+1}s^{(i+1)}_{i}
 \cdots s^{(i+1)}_{0}$
where the digits  $s^{(i+1)}_{0}=s_0$, \ldots, $s^{(i+1)}_{i}=s_{i}$ are integer digits of the output,
and
the factor $s^{(i+1)}_{i+3}s^{(i+1)}_{i+2}s^{(i+1)}_{i+1}$ is under consideration.
Note that for $j \ge i+3$, the coefficients $s^{(i+1)}_{j}$ are all equal to $0$.
Thus the next iteration of the algorithm gives
$\psi(s^{(i+1)}_{i+3}s^{(i+1)}_{i+2}s^{(i+1)}_{i+1})=
s^{(i+2)}_{i+3}s^{(i+2)}_{i+2}s^{(i+2)}_{i+1}$.

\begin{lem}
After every step $i$ of the algorithm, the coefficients satisfy:
\begin{itemize}
\item 
$s^{(i+1)}_{0}=s_0$, \ldots, $s^{(i+1)}_{i}=s_{i}$ belong to $\mathcal{A}$
\item
$s^{(i+1)}_{i+1}\in (-1,a)$ 
\item
$s^{(i+1)}_{i+2}\in (-1,0]$.
\end{itemize}
 \end{lem}
\begin{pf}
We will prove the statement by induction on the number of steps of the algorithm. 
The statement is valid for $i=0$ due to 
the assumption $|q|<1$.

By Transformation~\eqref{eq:alg_trans} we have
$\psi(s^{(i+1)}_{i+3}s^{(i+1)}_{i+2}s^{(i+1)}_{i+1})=s^{(i+2)}_{i+3}s^{(i+2)}_{i+2}s^{(i+2)}_{i+1}$, thus
\begin{align*}
s^{(i+2)}_{i+1}&= \lceil s^{(i+1)}_{i+1}\rceil  \in \Zset\cap[0,a]\,, \\
s^{(i+2)}_{i+2}&= s^{(i+1)}_{i+2} + a(\lceil s^{(i+1)}_{i+1}\rceil -s^{(i+1)}_{i+1}) \in (-1,a)\,, \\
s^{(i+2)}_{i+3} &= - (\lceil s^{(i+1)}_{i+1}\rceil -s^{(i+1)}_{i+1}) \in (-1,0]\,.
\end{align*}
\qed \end{pf}

Since the factor $s^{(i+2)}_{i+3}s^{(i+2)}_{i+2}s^{(i+2)}_{i+1}$ after step $i+1$ is
uniquely determined by the factor
$s^{(i+1)}_{i+3}s^{(i+1)}_{i+2}s^{(i+1)}_{i+1}$, and the coefficients $s^{(i+1)}_{i+1}$ and
$s^{(i+1)}_{i+2}$ are uniformly bounded, as a corollary we get
the following result.

\begin{prop}
Algorithm~\ref{algo:aexp_q} generates an $\alpha$-adic representation of $q$ which is
eventually periodic.
\end{prop}

\begin{exmp}\label{31}
We compute an $\alpha$-adic representation of the number $\frac{1}{2}$ in the case $\db=31$.
\begin{displaymath}
\begin{array}{r@{\ }r@{\ }r@{\ }r@{\ }r@{\ }r}
 & & & & & \phantom{-}\frac{1}{2} \\
 & & & -\frac{1}{2} & \phantom{-}\frac{3}{2} & \frac{1}{2}  \\
\hline
 & & & -\frac{1}{2} & \frac{3}{2} & 1 \\
 & & -\frac{1}{2} & \frac{3}{2} & \frac{1}{2} & \phantom{1} \\
\hline
 & & -\frac{1}{2} & 1 & 2 & 1 \\
-\frac{1}{2} & \phantom{-}\frac{3}{2} & \frac{1}{2} \\
\hline
-\frac{1}{2} & \phantom{-}\frac{3}{2} & 0 & 1 & 2 & 1
\end{array}
\end{displaymath}
Because the prefix $(-\frac{1}{2})(\frac{3}{2})$ which arises after step 3 is the same as
the one which arises after step 0, the same sequence of steps (with the same results) will follow from
now on. Therefore the $\alpha$-adic representation computed by the algorithm is 
$\aexp{\frac{1}{2}} = \lomega(012)1\fdot$.
It happens that, in this particular case, this is an $\alpha$-adic expansion of $\aexp{\frac{1}{2}}$.
\end{exmp}

\subsection{Normalization}

Unfortunately, Algorithm~\ref{algo:aexp_q} does not give directly an 
admissible $\alpha$-adic expansion in general.
In this section we discuss the normalization of such a non-admissible output.
Since the output word is a word on the canonical alphabet, its only possible non-admissible
factors are either of the type $a^nb$ with $n\geq 1$, $b\neq a$ or of the type $\lomega a$. 
The following result shows that the latter case will not appear.

\begin{prop}\label{prop:bounded_a}
The number of consecutive letters $a$'s in an output word of Algorithm~\ref{algo:aexp_q}
is bounded for all $q\in\Qset\cap(-1,1)$.
\end{prop}
\begin{pf}
We will prove the result by contradiction. Let us assume that from some step on,
say from step $i$, the output of the algorithm is composed only of letters $a$'s. 
This means that the output is of the form 
$\cdots\lceil V_4\rceil\lceil V_3\rceil\lceil V_2\rceil\lceil V_1\rceil v$, where
$v$ has length $i+1$, and for each $k \ge 1$, $\lceil V_k\rceil=a$.
We have
$V_1=s^{(i+1)}_{i+1}$, and $V_2=s^{(i+1)}_{i+2}+a(\lceil V_1\rceil -V_1)$.
Iterating twice the transformation $\psi$, we get
\begin{equation}\label{eq:Vk}
V_k = -(\lceil V_{k-2}\rceil - V_{k-2}) + a(\lceil V_{k-1}\rceil -V_{k-1})\quad \textrm{for}\quad k\geq 3.
\end{equation}

{}From Relation~\eqref{eq:Vk} and the fact that $V_k >a-1$ one get
\begin{equation}\label{eq:Vk_rec}
1 - \frac{1}{a} + \frac{1}{a}(\lceil V_{k-2}\rceil - V_{k-2}) < (\lceil V_{k-1}\rceil - V_{k-1})\,.
\end{equation}
Then iterating ~\eqref{eq:Vk_rec} we obtain an explicit estimate for $(\lceil V_k\rceil - V_k)$
\begin{align*}
(\lceil V_k\rceil - V_k) &> 1 - \frac{1}{a} + \frac{1}{a}(\lceil V_{k-1}\rceil - V_{k-1} ) \\
  &> 1 - \frac{1}{a} + \frac{1}{a}\Big( 1 - \frac{1}{a} + \frac{1}{a}(\lceil V_{k-2}\rceil - V_{k-2} )\Big) \\
  &= 1 - \frac{1}{a^2} + \frac{1}{a^2} (\lceil V_{k-2}\rceil - V_{k-2} )\\
  &> 1 - \frac{1}{a^3} + \frac{1}{a^3} (\lceil V_{k-3}\rceil - V_{k-3} )\\
  &> \cdots \\
  &> 1 - \frac{1}{a^{k-1}} + \frac{1}{a^{k-1}} (\lceil V_{1}\rceil - V_{1} )
\end{align*}
Since $s^{(i+1)}_{i+2} \in (-1,0]$ we can estimate 
$a-1 < V_2=s^{(i+1)}_{i+2}+a(\lceil V_1\rceil -V_1) \leq a(\lceil V_1\rceil -V_1)$,
which gives $1 - \frac{1}{a} <  (\lceil V_1\rceil - V_1)$. Therefore we have
\begin{equation}\label{eq:est_t1}
- (\lceil V_k\rceil - V_k) < -1 + \frac{1}{a^{k-1}} -\frac{1}{a^{k-1}} \Big(1 - \frac{1}{a}\Big) = \frac{1}{a^k} - 1\,.
\end{equation}
Finally, by inequality~\eqref{eq:est_t1} and the fact that $a-1 <V_k \le a$, we obtain a bound on $V_k$
\begin{equation}\label{eq:fin_est}
V_k = \underbrace{-(\lceil V_{k-2}\rceil - V_{k-2})}_{<\frac{1}{a^{k-2}} - 1} \
     \underbrace{+ a\lceil V_{k-1}\rceil}_{= a^2} \
     \underbrace{-aV_{k-1}}_{< a-a^2} < a - 1 + \frac{1}{a^{k-2}}\,.
\end{equation}

Suppose that we are computing an $\alpha$-adic expansion of a rational number $q$ with
denominator $p\in\Nset$. Find the smallest
$K$ such that $\frac{1}{p} > \frac{1}{a^{K-2}}$. Since any $V_k$ is a fraction with
denominator $p$, by~\eqref{eq:fin_est} we have $V_K = \frac{t}{p} < a - 1 + \frac{1}{a^{K-2}}$,
which implies $V_K < a-1$. This is in contradiction with the assumption that $a-1 < V_k$ for 
all $k\geq 1$.
\qed \end{pf}

\begin{prop}
Let $w$ be an output of Algorithm~\ref{algo:aexp_q} for a number $q\in\Qset\cap(-1,1)$
and let $\widehat{w}$ be the image of $w$ under the normalization function, $\nu_{\mathcal{A}}(w)=\widehat{w}$. 
Then $\widehat{w}$ is left eventually periodic with no fractional part.
\end{prop}
\begin{pf}
First of all, a number $\beta$ such that $\db = a1$ is a so-called \emph{confluent Pisot number} 
(cf.~\cite{frougny-tcs-106}). For these numbers,
it is known that the  normalization on the canonical alphabet does not produce a
carry to the right. This assures that $\widehat{w}$ will have no fractional part and that we
can perform normalization starting from the fractional point and then just read and write from right
to left. 

We have shown earlier that for a given rational number $q$ the number of consecutive
letters $a$'s in an output word $w$ is bounded, moreover the proof of 
Proposition~\ref{prop:bounded_a} gives us this upper bound. We give here a construction
of a right sequential transducer $\mathcal{T}$ performing the normalization of such a word $w$.

Define $\mathcal{A}_{\emptyset} \coloneq \mathcal{A}\setminus\{0\}$, and let $C$ be the bound on the 
number of consecutive letters $a$ in a word $w$.
Because the result of the normalization of non-admissible factors of $w$ depends on the
parity of the length of blocks of consecutive $a$'s, the transducer $\mathcal{T}$ has to 
count this parity. This is done by memorizing the actually processed forbidden factors;
the states of the transducer are labeled by these memorized factors.

Transducer $\mathcal{T}$ is constructed as follows
\begin{itemize}
\item The initial state is labeled by the empty word $\varepsilon$, and there is a loop
$\varepsilon \xrightarrow{0|0}\varepsilon$. 
\item There are states labeled by a single letter $h\in\mathcal{A}_{\emptyset}$ connected 
with the initial state by edges $\varepsilon \xrightarrow{h|\varepsilon}h$ and 
$h\xrightarrow{0|0h}\varepsilon$. These states are also connected one with each other
by edges $i \xrightarrow{j|i} j$ where $i,j\in\mathcal{A}_{\emptyset}$, $j\neq a$.
Finally there is a loop $h\xrightarrow{h|h}h$ on each state $h\in\mathcal{A}_{\emptyset}$, $h\neq a$.
\item For each $h\in\mathcal{A}_{\emptyset}$ there is a chain of consecutive states $a^kh$, where $k=1,\ldots,C-1$,
linked by edges $a^kh \xrightarrow{a|\varepsilon} a^{k+1}h$. Moreover, there are edges $a^kh \xrightarrow{i|u}i+1$
where $u = (0a)^m0(h-1)$ for $k=2m+1$ and $u = (0a)^m0(a-1)h$
for $k=2m+2$.
\end{itemize}
The edges $a^kh \xrightarrow{a|\varepsilon} a^{l+1}h$ are these which count the number of consecutive letters
$a$ in a forbidden factor, whereas the edges $a^kh \xrightarrow{i|u}i+1$ are these which, depending on the
parity of the length $k$ of a
run $a^k$, replace a forbidden factor by its normalized equivalent.

One can easily check that the transducer is input deterministic, and thus right sequential.
Clearly the output word is admissible.
Since the image by a sequential function of an eventually periodic
word is eventually periodic (see~\cite{eilenberg}),
the image $\hat{w}$ is
eventually periodic. 
\qed \end{pf}

The following is just a rephrasing.
\begin{thm}
Let $\beta$ be a quadratic Pisot unit. Any rational number $q \in\Qset\cap(-1,1)$ has an eventually
periodic $\alpha$-adic expansion with no fractional part.
\end{thm}

Remark that there exist rational numbers larger than $1$ such that the $\alpha$-adic expansion
has no fractional part.
We have shown in Example~\ref{31} that for 
$\db=31$, $\aexp{\frac{1}{2}} = \lomega(012)1\fdot$. Thus $\aexp{\frac{3}{2}} = \lomega(012)2\fdot$
has no fractional part.

\section{Conclusion}
Let us stress out that the analogue of Propositions \ref{prop:z_plus2} and \ref{quad-} has
 been proved by Sadahiro for the case that
$\beta$ is a cubic Pisot unit with complex conjugates satisfying Property (F).
The extension of these results to other Pisot units satisfying Property (F) is an open problem.

\begin{ack}

The authors are grateful to Shigeki Akiyama and Christoph Bandt for stimulating discussions.

\end{ack}


\begin{thebibliography}{10}

\bibitem{akiyama-dm-2}
S.~Akiyama.
\newblock {\em Self affine tiling and {P}isot numeration system}.
\newblock In 'Number theory and its applications (Kyoto, 1997)', K.~Gy{\"o}ry
  and S.~Kanemitsu, (eds.), Dev. Math. {\bf 2}, Kluwer Acad. Publ.
  (1999), 7--17.

\bibitem{akiyama-00}
S.~Akiyama.
\newblock {\em Cubic {P}isot units with finite beta expansions}.
\newblock In 'Algebraic number theory and Diophantine analysis (Graz, 1998)',
  de Gruyter (2000), 11--26.

\bibitem{akiyama-sadahiro-amiuo-6}
S.~Akiyama and T.~Sadahiro.
\newblock {\em A self-similar tiling generated by the minimal {P}isot number}.
\newblock In 'Proceedings of the 13th Czech and Slovak International Conference
  on Number Theory (Ostravice, 1997)', Acta Math. Inform. Univ. Ostraviensis
{\bf 6}, 9--26, (1998).

\bibitem{ambroz-words-05}
P.~Ambro\v{z}.
\newblock {\em On the tau-adic expansions of real numbers}.
\newblock In 'Words 2005, 5{$^{th}$} International Conference on Words, actes',
  S.~Brlek and C.~Reutenauer, (eds.), Publications du LaCIM {\bf 36},
  UQ\`{A}M (2005), 79--89.

\bibitem{berthe-siegel-arxiv02}
V.~Berth{\'e} and A.~Siegel.
\newblock {\em Purely periodic beta-expansions in the {P}isot non-unit case}.
\newblock  Rapport de recherche LIRMM 04025, Arxiv math. DS/0407282, (2002).

\bibitem{bertrand-asp-285}
A.~Bertrand.
\newblock {\em D\'eveloppements en base de {P}isot et r\'epartition modulo {$1$}}.
\newblock C. R. Acad. Sci. Paris {\bf 285} (1977), 419--421.

\bibitem{blanchard-tcs-65}
F.~Blanchard.
\newblock {\em $\beta$-expansions and symbolic dynamics}.
\newblock Theoret. Comput. Sci. {\bf 65} (1989), 131--141.

\bibitem{eilenberg}
S.~Eilenberg.
\newblock {\em Automata, Languages, and Machines}
Vol. A.
\newblock Academic Press, 1974.

\bibitem{frougny-tcs-106}
Ch.~Frougny.
\newblock {\em Confluent linear numeration systems}.
\newblock Theoret. Comput. Sci. {\bf 106} (1992), 183--219.

\bibitem{frougny-mst-25}
Ch.~Frougny.
\newblock {\em Representations of numbers and finite automata}.
\newblock Math. Systems Theory {\bf 25} (1992), 37--60.

\bibitem{frougny02} Ch.~Frougny.
\newblock {\em On-line odometers for two-sided symbolic
dynamical systems}.
\newblock Proceedings of DLT 2002, Lecture Notes in Computer Science {\bf 2450} (2002), 405--416.

\bibitem{frougny-solomyak-etds-12}
Ch.~Frougny and B.~Solomyak.
\newblock {\em Finite beta-expansions}.
\newblock Ergod. Th. and Dynam. Sys. {\bf 12} (1992), 713--723.

\bibitem{grabner-liardet-tichy-aa-80}
P.~Grabner, P.~Liardet, and R.~Tichy.
\newblock {\em Odometers and systems of numeration}.
\newblock Acta Arith. {\bf 80} (1995), 103--123.

\bibitem{hollander-thesis}
M.~Hollander.
\newblock {\em Linear numeration systems, finite beta-expansions, and discrete
  spectrum of substitution dynamical systems}.
\newblock PhD thesis, Washington University, (1996).

\bibitem{ito-rao-pams-04}
S.~Ito and H.~Rao.
\newblock {\em Purely periodic {$\beta$}-expansions with {P}isot unit base}.
\newblock Proc. of Amer. Math. Soc. {\bf 133} (2004), 953--964.

\bibitem{katai}
I.~K\'atai.
\newblock {\em Number systems in imaginary quadratic fields}.
\newblock Ann. Univ. Sci. Budapest Sect. Comput. {\bf 14} (1994), 91--103.

\bibitem{lothaire2}
M.~Lothaire.
\newblock {\em Algebraic {C}ombinatorics on {W}ords}, volume~90 of 
  {\em Encyclopedia of Mathematics and its Applications}.
\newblock Cambridge University Press, Cambridge, (2002).

\bibitem{messaoudi-aa-95}
A.~Messaoudi.
\newblock {\em Fronti\`ere du fractal de {R}auzy et syst\`eme de num\'eration complexe}.
\newblock Acta Arith. {\bf 95} (2000), 195--224.

\bibitem{parry-amash-11}
W.~Parry.
\newblock {\em On the {$\beta$}-expansions of real numbers}.
\newblock Acta Math. Acad. Sci. Hungar. {\bf 11} (1960), 401--416.

\bibitem{rauzy-bsmf-110}
G.~Rauzy.
\newblock {\em Nombres alg\'ebriques et substitutions}.
\newblock Bull. Soc. Math. France {\bf 110} (1982), 147--178.

\bibitem{renyi-amash-8}
A.~R{\'e}nyi.
\newblock {\em Representations for real numbers and their ergodic properties}.
\newblock Acta Math. Acad. Sci. Hungar {\bf 8} (1957), 477--493.

\bibitem{sadahiro-pp}
T.~Sadahiro.
\newblock {\em Multiply covered points of dual Pisot tilings}.
\newblock Preprint, (2005).

\bibitem{schmidt-blms-12}
K.~Schmidt.
\newblock {\em On periodic expansions of {P}isot numbers and {S}alem numbers}.
\newblock Bull. London Math. Soc. {\bf 12} (1980), 269--278.

\bibitem{schmidt-mm-129}
K.~Schmidt.
\newblock {\em Algebraic coding of expansive group automorphisms and two-sided beta-shifts}.
\newblock Monatsh. Math. {\bf 129} (2000), 37--61.

\bibitem{sidorov-vershik-mm-126}
N.~Sidorov and A.~Vershik.
\newblock {\em Ergodic properties of the {E}rd{\H o}s measure, the entropy of the golden shift, 
  and related problems}.
\newblock Monatsh. Math. {\bf 126} (1998), 215--261.

\bibitem{vershik-iml-92}
A.~M. Vershik.
\newblock {\em The fibadic expansions of real numbers and adic transformations}.
\newblock Prep. Report Inst. Mittag-Leffler {\bf } (1991/1992), 1--9.

\end{thebibliography}
\end{document}